\documentclass[11pt]{article}
\usepackage[paper=a4paper,left=20mm,width=170mm,top=20mm,bottom=25mm]{geometry}
\usepackage[english]{babel}
\usepackage{amsmath,amsthm}
\usepackage{amsfonts}
\usepackage{enumitem}
\usepackage{wasysym}
\usepackage{lipsum}
\usepackage{framed}
\usepackage{color}
\usepackage{tikz}
\usepackage{tikz-qtree}
\usepackage{url}
\usetikzlibrary{trees}

\theoremstyle{definition}

\theoremstyle{remark}

\numberwithin{equation}{section}
\DeclareMathOperator{\sgn}{sgn}
\theoremstyle{acknowledgment}
\newtheorem*{acknowledgment}{Acknowledgment}

\begin{document}

\title{\bfseries Integral of radical trigonometric functions revisited}
\author{\normalsize N. Karjanto$^{1}$ and B. Yermukanova$^{2}$\\
{\small $^{1}$Department of Mathematics, University College, Sungkyunkwan University}\\
{\small Natural Science Campus, Suwon, Republic of Korea}\\
{\small $^{2}$Department of Economics, School of Humanities and Social Sciences}\\
{\small Nazarbayev University, Astana, Kazakhstan}}
\date{\scriptsize Last updated \today}%
\date{\scriptsize }%
\maketitle

\begin{abstract}
This article revisits an integral of radical trigonometric functions.
It presents several methods of integration where the integrand takes the form $\sqrt{1 \pm \sin x}$ or $\sqrt{1 \pm \cos x}$.
The integral has applications in Calculus where as the length of cardioid represented in polar coordinates.\\

Keywords: Techniques of integration; radical trigonometric functions; cardioid.
\end{abstract}

\section{Introduction}
This article revisits an integral where the integrand takes the form of radical trigonometric functions.
A general form of radical trigonometric integrands in the context of this article refers to $\sqrt{a \pm b \sin x}$ or $\sqrt{a \pm b \cos x}$, for $a, b > 0$. 
The integral of these functions is expressed in terms of elliptic integral and are available in mathematical handbooks and tables of integrals. 
For example, the latter integral is given in Section 2.5 (see 2.576) of a famous mathematical handbook by Gradstyen and Ryzhik~\cite{gradshteyn}.
For a particular case of $a = b$, after removing the constant factor, the integrand reduces to radical trigonometric functions $\sqrt{1 \pm \sin x}$ or $\sqrt{1 \pm \cos x}$. Interestingly, it seems that explicit expressions for the integral of these functions have not been specifically listed in any tables of integrals and handbooks, including, but not limited to,~\cite{abramowitz,bronshtein,gradshteyn,table,mathar,polyanin,spiegel}. 
The focus of this article is to consider the special case when $a = 1 = b$, where several techniques of integration are discussed in a more detail. 
To the best of our knowledge, this is the first time when such a compilation for particular integrands is presented. 

\begin{figure}[h]
\begin{center}
\includegraphics[width=0.99\textwidth]{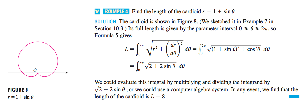}
\end{center}
\caption{An example from a textbook on calculating the length of cardioid $r = (1 + \sin \theta)$ where the calculation details are omitted.} \label{example}
\end{figure}

The motivation of this article springs from an encounter from one of the coauthors' in teaching Calculus~2 course during the Spring 2014 semester in 
Nazarbayev University, Astana, Kazakhstan. 
Particularly, the content of this article is related to the topic on the integral calculus of polar curves, and one of the examples is calculating the length of a cardioid. We adopt the Calculus textbook written by Anton, Bivens and Davis~\cite{anton} where the polar curves are discussed in Section 10.3. Another recommended textbook reading for this course is the one written by Stewart~\cite{stewart12}. An example [Example 4, Section 10.4, page 692] from the latter textbook mentions that finding the length of cardioid $r = 1 + \sin \theta$ can be evaluated by multiplying both the numerator and the denominator of the integrand by $\sqrt{2 - 2 \sin \theta}$ or alternatively, using the Computer Algebra System (CAS). Yet, evaluating this integral by hand is apparently not so obvious to many students since they have to do further manipulation on the obtained expression. The screenshot of the example from Stewart's textbook has been excerpted and displayed in Figure~\ref{example}.

After rationalizing the numerator and implementing the Pythagorean trigonometric identity, the numerator simplifies to $\sqrt{\cos^2\theta} = |\cos \theta|$, but it has to be in absolute value form, instead of simply $\cos \theta$. 
This is a common mistake found among students since they may forget or tend to ignore the absolute value sign.
From an instructor's perspective, it is imperative to remind the students to be aware of this fact. Referring to Bloom's taxonomy of learning domains~\cite{bloom}, the educational activity of this learning process is the {\sl cognitive domain}. The process covers {\sl knowledge}, {\sl comprehension} and {\sl application}. In this example, students possess the knowledge that any value of the square root must be non-negative and an absolute value of any quantity is always non-negative too. A comprehension of these facts is essential to conclude (application of knowledge) that the square root of a quantity squared is indeed equal to the absolute value of that quantity. 

Referring to the revised Bloom taxonomy~\cite{anderson}, a connection between learning activities and learning objectives can further be established.
The knowledge dimension covers the {\sl factual} and {\sl conceptual} aspects. 
In this context, students must know the definition of an absolute value and be able to make an interrelationship between the property of a square root and the absolute value.
The cognitive dimension includes {\sl remember, understand, apply} aspects. Possessing the knowledge of absolute value, it is crucial to investigate whether the students can retrieve this knowledge from their memory, whether they understand why absolute value has to be non-negative and whether they are able to simplify and conclude that $\sqrt{\cos^2 \theta} = |\cos \theta|$.

Another educational aspect of the integral involving radical trigonometric functions is related to the {\sl synthesis} skill of {\sl cognitive domain} in Bloom's taxonomy. In the revised Bloom's taxonomy, the educational content involves {\sl factual} and {\sl conceptual} aspects of the {\sl knowledge} dimension, where students attempt to make interrelationships among the basic elements of trigonometric functions. The {\sl cognitive process} dimension covers {\sl remember, understand, apply} and {\sl analyze} aspects. When this example is posed to the classroom for the students to work on, it turns out that some excellent students come up with different techniques by manipulating the integrand expression. This shows that different students approach the problem distinctly, they attempt to integrate with the method which is most convenient to them. For instance, for some students, the technique of trigonometric substitution is a more comfortable approach, others implement a variable shift method to solve the problem successfully. Thus, there are several ways by which students can approach the problem. 

For many Calculus instructors, however, the interest in integration techniques has waned. 
With the introduction of CAS, many of them now give only cursory attention to such techniques.
Nevertheless, the methods of integration covered in this article are still interesting from educational perspective.
They provide a valuable pedagogical tool to assist and improve the students' learning skills, which are beneficial to both the instructors and the students themselves alike.
In particular, by introducing several methods during class sessions, the techniques covered in this article become useful in the sense that it does not only expose the students to various techniques of integration but also makes them review and strengthen their knowledge of trigonometry and trigonometric functions.	
As can be observed later, this article recalls some important properties of trigonometric functions of sine, cosine and tangent as well as a significant application of trigonometric substitution in solving particular types of integration.

This article is organized as follows.
The following section covers the integral of radical sine function $\sqrt{1 \pm \sin x}$.
Section~\ref{cosine} briefly covers the integral of radical cosine function $\sqrt{1 \pm \cos x}$.
Several techniques of integration are covered and more detailed derivations are discussed in Section~\ref{sine}, 
including rationalizing numerator, combining trigonometric identities, twice trigonometric substitutions and variable shift methods.
All of these methods require some variations of integrating absolute value function, which will be presented accordingly in the corresponding subsections. Section~\ref{cardioid} presents an application where the integrals of radical sine and cosine functions appear, particularly in calculating the length of a cardioid. The final section draws conclusions and provides remark to our discussion.

\section{Integral of radical sine function} \label{sine}

This section deals with the integral of a radical sine function where the integrand takes the form $\sqrt{1 \pm \sin x}$.
There are a number of methods to obtain the result, and four techniques are covered in this section.
The first method is by rationalizing the numerator. From here, one may depart either to use the Pythagorean identity or to employ a trigonometric substitution.
The second method is by combining several trigonometric identities. We observe that double-angle formula and the identity relating $\sin x$ and $\tan x/2$ 
follow different paths of calculation and yet arrive at identical expression. 
The third technique is by implementing trigonometric substitutions two times, mainly using tangent function.
Finally, the fourth technique is conducted by shifting the variable by $\pi/2$. Two options can be developed from this path, where
both of them alter the integral from radical sine function into radical cosine function.
The four methods covered in this section are summarized in the following tree diagram.
\begin{center}
\begin{tikzpicture}[level distance=6cm,sibling distance=0.5cm,scale=.8]
\tikzset{edge from parent/.style= 
          {thick, draw, edge from parent fork right}, 
          every tree node/.style = {draw = blue!80, shape = rectangle, rounded corners =.8ex, fill = cyan!30, 
          minimum width=5cm,text width=5cm, align=center},grow'=right}
\Tree 
    [. {Integration techniques}
        [.{Rationalizing numerator}
            [.{Pythagorean identity} ]
            [.{Trigonometric substitution} ]
        ]
        [.{Combining identities}
            [.{Double-angle formula} ]
            [.{$\sin x$ and $\tan (x/2)$} ]
        ] 
        [.{Twice trigonometric\\ substitutions} ]
        [.{Variable shift}
            [.{$x = y - \pi/2$} ]
            [.{$x = \pi/2 - y$} ]
        ]
    ]
\end{tikzpicture}
\end{center}

\subsection{Rationalizing numerator} \label{ratnum}
The following integral will be used in this subsection.
Let $f$ be a function which has at most one root on each interval on which it is defined, and $F$ an antiderivative of $f$, i.e. $F'(x) = f(x)$, then
\begin{equation}
\int \frac{\left| f(x)\right|}{\sqrt{F(x)}} \,dx = -2 \sgn[f(x)] \, \sqrt{F(x)} + C  \label{av1}
\end{equation}
where $\sgn(x)$ is the sign function, which takes the values $-1, 0$ or $1$ when $x$ is negative, zero or positive, respectively.

\subsubsection*{Pythagorean identity}
Let $I$ be an indefinite integral of the radical sine function 
${\displaystyle  I = \int \sqrt{1 \pm \sin x}\, dx}$,
then rationalizing the numerator by multiplying both the numerator and the denominator with $\sqrt{1 \mp \sin x}$, applying the Pythagorean trigonometric identity and utilizing the definition of the absolute value, it yields:
\begin{eqnarray*}
  I &=& \int \sqrt{1 \pm \sin x} \cdot \frac{\sqrt{1 \mp \sin x}}{\sqrt{1 \mp \sin x}} \, dx  
     =  \int \frac{\sqrt{1 - \sin^2 x}}{\sqrt{1 \mp \sin x}} \, dx 
     =  \int \frac{\sqrt{\cos^2 x}}{\sqrt{1 \mp \sin x}} \, dx  \\
    &=& \int \frac{|\cos x|}{\sqrt{1 \mp \sin x}} \, dx 
     = -2 \sgn(\cos x) \sqrt{1 \mp \sin x} + C
\end{eqnarray*}
where the last expression is readily obtained by implementing~\eqref{av1}.

\subsubsection*{Trigonometric substitution}

A similar solution can also be obtained using the trigonometric substitution of $u = \sin x$.
Differentiating with respect to $u$, we get $dx = du/\cos x = du/(\pm \sqrt{1 - u^2})$, where the positive and negative signs are related to the sign of $\cos x$. Thus for $u = \sin x \neq \pm 1$ 
\begin{eqnarray*}
  I &=& \int \frac{\sqrt{1 \pm u}\, du}{\pm \sqrt{1 - u^2}} 
     =  \int \frac{du}{\pm \sqrt{1 \mp u}} 
     =  \mp 2 \sqrt{1 \mp u} + C \\
    &=& \mp 2 \sqrt{1 \mp \sin x} + C 
     =  -2 \sgn(\cos x) \sqrt{1 \mp \sin x} + C.
\end{eqnarray*}

\subsection{Combining identities}
A general, explicit form of an integral involving an absolute value of a function will be used in this section.
Let $f$ be a function which has at most one root on each interval on which it is defined, and $F$ an antiderivative of $f$ that is zero at each root of $f$ (such an antiderivative exists if and only if the condition on $f$ is satisfied), then
\begin{equation}
\int \left| f(x)\right| \,dx = \sgn[f(x)] \, F(x) + C, \label{av2}
\end{equation}
where sgn$(x)$ is the sign function defined previously.

\subsubsection*{Double-angle formula}
We manipulate the integrand by combining the Pythagorean trigonometric identity and the double-angle formula.
Using the Pythagorean trigonometric identity, writing $1 = \cos^2 (x/2) + \sin^2 (x/2)$ and using the double-angle formula for $\sin x$: $\sin x = 2 \sin (x/2) \cos(x/2)$, the integral of the radical sine becomes
\begin{eqnarray*}
  I &=& \int \sqrt{\cos^2\frac{x}{2} \pm 2 \cos \frac{x}{2} \sin \frac{x}{2} + \sin^2\frac{x}{2}} \, dx 
     =  \int \sqrt{\left(\cos \frac{x}{2} \pm \sin \frac{x}{2} \right)^2}\, dx \\
    &=& \int \left|\cos \frac {x}{2} \pm \sin \frac{x}{2} \right| \, dx 
     =  2 \sgn \left(\cos \frac{x}{2} \pm \sin \frac{x}{2} \right) \left(\sin \frac{x}{2} \mp \cos \frac{x}{2} \right) + C
\end{eqnarray*}
where the last expression is quickly obtained after implementing~\eqref{av2}.

\subsubsection*{Identity relating $\sin x$ and $\tan (x/2)$}
A similar result will also be obtained if one employs another trigonometric identity that relates $\sin x$ and $\tan (x/2)$.
Using the double-angle formula for $\sin x$ at the numerator and the Pythagorean trigonometric identity at the denominator, dividing both sides by $\cos^{2} (x/2)$, we obtain
\begin{equation*}
  \sin x = \frac{2 \sin (x/2) \cos(x/2)}{\cos^2 (x/2) + \sin^2 (x/2)} = \frac{\frac{2 \sin (x/2) \cos (x/2)}{\cos^2(x/2)}}{1 + \frac{\sin^2 (x/2)}{\cos^2 (x/2)}} = \frac{2 \tan (x/2)}{1 + \tan^2 (x/2)}.
\end{equation*}
Thus, the integral of the radical sine function $I$ turns to
\begin{eqnarray*}
  I &=& \int \sqrt{1 \pm \frac{2 \tan x/2}{1 + \tan^2 x/2}} \, dx 
     =  \int \sqrt{\frac{1 + \tan^2 x/2 \pm 2 \tan x/2}{1 + \tan^2 x/2}} \, dx \\
    &=& \int \sqrt{\frac{(1 \pm \tan x/2)^2}{\sec^2 x/2}} \, dx 
     =  \int \left|\frac{1 \pm \tan x/2}{\sec x/2} \right| \, dx \\
    &=& \int \left|\cos \frac {x}{2} \pm \sin \frac{x}{2} \right| \, dx 
     =  2 \sgn \left(\cos \frac{x}{2} \pm \sin \frac{x}{2} \right) \left(\sin \frac{x}{2} \mp \cos \frac{x}{2} \right) + C.
\end{eqnarray*}

\subsection{Twice trigonometric substitutions}
A similar expression of the solution as that of the previous section can also be obtained by the trigonometric substitution $u = \tan x/2$.
This implies $dx = 2 \, du/(1 + u^2)$ and writing $\sin x = 2 \sin (x/2) \cos (x/2)$ the integral of the radical sine function becomes
\begin{eqnarray*}
  I &=& \int \sqrt{1 \pm \frac{2u}{1 + u^2}} \, \frac{2du}{1 + u^2} 
     =  2 \int \frac{|1 \pm u| \, du}{(1 + u^2)^{3/2}} \\
    &=& 2 \sgn(1 \pm u) \left(\int \frac{du}{(1 + u^2)^{3/2}} \pm \int  \frac{u \, du}{(1 + u^2)^{3/2}} \right).
\end{eqnarray*}
Employ another trigonometric substitution $u = \tan y$ and $v = 1 + u^2$ for the first and the second integrals, respectively.
Thus,
\begin{eqnarray*}
  I &=& 2 \sgn(1 \pm u) \left(\int \frac{\sec^2 y\, dy}{(1 + \tan^2 y)^{3/2}} \pm \frac{1}{2} \int \frac{dv}{v^{3/2}} \right) 
     =  2 \sgn(1 \pm u) \left(\int \frac{\sec^2 y\, dy}{\sec^3 y} \mp v^{-1/2} \right) \\
    &=& 2 \sgn(1 \pm u) \left(\int \frac{1}{\sec y} \, dy \mp \frac{1}{\sqrt{v}} \right) 
     =  2 \sgn(1 \pm u) \left(\int \cos y \, dy \mp \frac{1}{\sqrt{1 + u^2}} \right) \\
    &=& 2 \sgn(1 \pm u) \left(\sin y \mp \frac{1}{\sqrt{1 + u^2}} \right) + C 
     =  2 \sgn(1 \pm u) \sgn\left(\frac{1}{\sqrt{1 + u^2}} \right) \left(\frac{u \mp 1}{\sqrt{1 + u^2}} \right) + C\\
    &=& 2 \sgn\left(1 \pm \tan\frac{x}{2} \right) \sgn \left(\cos \frac{x}{2} \right) \; \cos \frac{x}{2} \left(\tan \frac{x}{2} \mp 1 \right) + C \\ 
    &=& 2 \sgn \left(\cos \frac{x}{2} \pm \sin \frac{x}{2} \right) \left(\sin \frac{x}{2} \mp \cos \frac{x}{2} \right) + C.
\end{eqnarray*}

\subsection{Variable shift}
The following integrals of the absolute value of trigonometric functions $\sin{\alpha x}$ and $\cos{\alpha x}$, $\alpha \neq 0$, will be used in this subsection, where $\lfloor x \rfloor$ denotes the floor function:
\begin{eqnarray}
\int \left|\sin{\alpha x} \right|\,dx &=& {2 \over \alpha} \left\lfloor \frac{\alpha x}{\pi} \right\rfloor - {1 \over \alpha} \cos{\left(\alpha x - \left\lfloor \frac{\alpha x}{\pi} \right\rfloor \pi \right)} + C \label{av3}\\
\int \left|\cos{\alpha x} \right|\,dx &=& {2 \over \alpha} \left\lfloor \frac{\alpha x}{\pi} + \frac12 \right\rfloor + {1 \over \alpha} \sin{\left(\alpha x - \left\lfloor \frac{\alpha x}{\pi} + \frac12 \right\rfloor \pi \right)} + C. \label{av4}
\end{eqnarray}

\subsubsection*{Variable shift $x = y - \pi/2$}

Applying this variable shift, the integral $I$ becomes
\begin{eqnarray*}
  I &=& \int \sqrt{1 \pm \sin (y - \pi/2)} \, dy 
     =  \int \sqrt{1 \mp \cos y} \, dy \\
    &=&
\left\{
  \begin{array}{ll}
    {\displaystyle \int} \sqrt{2 \sin^2 y/2} \, dy, & \hbox{for $-$ sign ($+$ sign original $I$)} \\
    {\displaystyle \int} \sqrt{2 \cos^2 y/2} \, dy, & \hbox{for $+$ sign ($-$ sign original $I$)}
  \end{array}
\right.\\
&=&
\left\{
  \begin{array}{ll}
  \sqrt{2} {\displaystyle \int}  \left| \sin y/2 \right| \, dy, & \hbox{for $-$ sign ($+$ sign original $I$)} \\
  \sqrt{2} {\displaystyle \int}  \left| \cos y/2 \right| \, dy, & \hbox{for $+$ sign ($-$ sign original $I$)}
  \end{array}
\right.\\
&=&
\left\{
  \begin{array}{rl}
  -2 \sqrt{2} \sgn(\sin y/2) \cos(y/2) + C, & \hbox{for $-$ sign ($+$ sign original $I$)} \\
   2 \sqrt{2} \sgn(\cos y/2) \sin(y/2) + C, & \hbox{for $+$ sign ($-$ sign original $I$)}
  \end{array}
\right.\\
&=&
\left\{
  \begin{array}{rl}
  -2 \sqrt{2} \sgn\left[\sin \left(\frac{x}{2} + \frac{\pi}{4} \right) \right] \cos \left(\frac{x}{2} - \frac{\pi}{4} \right) + C, & \hbox{for $-$ sign ($+$ sign original $I$)} \\
   2 \sqrt{2} \sgn\left[\cos \left(\frac{x}{2} + \frac{\pi}{4} \right) \right] \sin \left(\frac{x}{2} + \frac{\pi}{4} \right) + C, & \hbox{for $+$ sign ($-$ sign original $I$)}
  \end{array}
\right.\\
    &=& -2 \sqrt{2} \sgn\left[\sin \left(\frac{x}{2} \pm \frac{\pi}{4} \right) \right] \cos \left(\frac{x}{2} \pm \frac{\pi}{4} \right) + C
\end{eqnarray*}
where the last three expressions are readily obtained by implementing~\eqref{av2}, returning back the original variable and combining results corresponding to the positive and negative signs into a single expression, respectively.
Alternatively, implementing~\eqref{av3}, we obtain the integral for $\sqrt{1 + \sin x}$:
\begin{eqnarray*}
I_1 &=& 4 \sqrt{2} \left\lfloor \frac{y}{2\pi} \right\rfloor - 2 \sqrt{2} \cos{\left({y \over 2} - \left\lfloor \frac{y}{2\pi} \right\rfloor \pi \right)} + C \\
    &=& 4 \sqrt{2} \left\lfloor \frac{x}{2\pi} + \frac14 \right\rfloor - 2 \sqrt{2} \cos{\left({x \over 2} + \frac{\pi}{4} - \left\lfloor \frac{x}{2\pi} + \frac{1}{4} \right\rfloor \pi \right)} + C.
\end{eqnarray*}
Implementing~\eqref{av4}, we obtain the integral for $\sqrt{1 - \sin x}$:
\begin{eqnarray*}
I_2 &=& 4 \sqrt{2} \left\lfloor \frac{y}{2\pi} + \frac12 \right\rfloor + 2 \sqrt{2} \sin{\left({y \over 2} - \left\lfloor \frac{y}{2\pi} + \frac12 \right\rfloor \pi \right)} + C \\
    &=& 4 \sqrt{2} \left\lfloor \frac{x}{2\pi} + \frac34 \right\rfloor + 2 \sqrt{2} \sin{\left({x \over 2} + \frac{\pi}{4} - \left\lfloor \frac{x}{2\pi} + \frac{3}{4} \right\rfloor \pi \right)} + C
\end{eqnarray*}
where subscripts 1 and 2 correspond to the positive and negative signs in the original integral $I$, respectively.

\subsubsection*{Variable shift $x = \pi/2 - y$}

Applying this variable shift, the integral $I$ becomes
\begin{eqnarray*}
  I &=& -\int \sqrt{1 \pm \sin (\pi/2 - y)} \, dy 
     =  -\int \sqrt{1 \pm \cos y} \, dy \\
&=&
\left\{
  \begin{array}{rl}
- \sqrt{2} {\displaystyle \int}  \left| \cos y/2 \right| \, dy, & \hbox{for $+$ sign} \\
- \sqrt{2} {\displaystyle \int}  \left| \sin y/2 \right| \, dy, & \hbox{for $-$ sign}
  \end{array}
\right.\\
&=&
\left\{
  \begin{array}{rl}
-2 \sqrt{2} \sgn(\cos y/2) \sin(y/2) + C, & \hbox{for $+$ sign} \\
 2 \sqrt{2} \sgn(\sin y/2) \cos(y/2) + C, & \hbox{for $-$ sign}
  \end{array}
\right.\\
&=&
\left\{
  \begin{array}{rl}
 2 \sqrt{2} \sgn\left[\cos \left(\frac{x}{2} - \frac{\pi}{4} \right) \right] \sin \left(\frac{x}{2} - \frac{\pi}{4} \right) + C, & \hbox{for $+$ sign} \\
-2 \sqrt{2} \sgn\left[\sin \left(\frac{x}{2} - \frac{\pi}{4} \right) \right] \cos \left(\frac{x}{2} - \frac{\pi}{4} \right) + C, & \hbox{for $-$ sign}
  \end{array}
\right.\\
    &=&  2 \sqrt{2} \sgn\left[\cos \left(\frac{x}{2} \mp \frac{\pi}{4} \right) \right] \sin \left(\frac{x}{2} \mp \frac{\pi}{4} \right) + C
  \end{eqnarray*}
where the last three expressions are readily obtained by implementing~\eqref{av2}, returning back the original variable and combining two results into a single expression, respectively.
Alternatively, implementing~\eqref{av4}, we obtain the integral for $\sqrt{1 + \sin x}$:
\begin{eqnarray*}
I_1 &=& -4 \sqrt{2}\left\lfloor \frac{y}{2\pi} + \frac12 \right\rfloor - 2 \sqrt{2}  \sin{\left(\frac{y}{2} - \left\lfloor \frac{y}{2\pi} + \frac12  \right\rfloor \pi \right)} + C \\
    &=& -4 \sqrt{2}\left\lfloor \frac34 -  \frac{x}{2\pi} \right\rfloor + 2 \sqrt{2}  \sin{\left(\frac{x}{2} - \frac{\pi}{4} + \left\lfloor \frac34 - \frac{x}{2\pi}  \right\rfloor \pi \right)} + C \\
    &=&  4 \sqrt{2}\left\lceil \frac{x}{2\pi} - \frac34 \right\rceil + 2 \sqrt{2}  \sin{\left(\frac{x}{2} - \frac{\pi}{4} - \left\lceil \frac{x}{2\pi} - \frac34  \right\rceil \pi \right)} + C
\end{eqnarray*}
where $\lceil x \rceil$ is the ceiling function and the relationship between the floor and the ceiling functions are utilized to obtain the last expression, i.e. $\lfloor x \rfloor + \lceil -x \rceil = 0$. Implementing~\eqref{av3}, we obtain the integral for $\sqrt{1 - \sin x}$:
\begin{eqnarray*}
I_2 &=& -4 \sqrt{2}\left\lfloor \frac{y}{2\pi} \right\rfloor + 2 \sqrt{2}  \cos {\left(\frac{y}{2} - \left\lfloor \frac{y}{2\pi}  \right\rfloor \pi \right)} + C \\
    &=& -4 \sqrt{2}\left\lfloor \frac14 -  \frac{x}{2\pi} \right\rfloor + 2 \sqrt{2}  \cos{\left(\frac{x}{2} - \frac{\pi}{4} + \left\lfloor \frac14 - \frac{x}{2\pi}  \right\rfloor \pi \right)} + C \\
    &=&  4 \sqrt{2}\left\lceil \frac{x}{2\pi} - \frac14 \right\rceil + 2 \sqrt{2}  \cos{\left(\frac{x}{2} - \frac{\pi}{4} - \left\lceil \frac{x}{2\pi} - \frac14  \right\rceil \pi \right)} + C
\end{eqnarray*}
where the subscripts 1 and 2 correspond to the positive and negative signs in the expressions of $I$, respectively.

\section{Integral of radical cosine function} \label{cosine}

This section compiles a number of techniques to integrate the radical cosine function in the form $\sqrt{1 \pm \cos x}$.
Let $J$ be an indefinite integral of radical cosine function
${\displaystyle J = \int \sqrt{1 \pm \cos x} \, dx.}$
Since the derivations are similar to the ones in Section~\ref{sine}, only the final results will be presented. 
Employing the variable shift method either by $x = \pi/2 - y$ or $x = y - \pi/2$ will alter the cosine function into the sine function and vice versa. 
Thanks to this redundancy, the coverage of this technique will be omitted in this section.
The integration techniques presented in this section basically can also be summarized with a similar tree diagram presented in Section~\ref{sine}.

\subsection{Rationalizing numerator}
Implementing two techniques of rationalizing numerator and by trigonometric substitution $u = \cos x$, 
we obtain a similar result to the one in the previous section:
\begin{equation*}
J = -2 \sgn(\sin x) \sqrt{1 \mp \cos x} + C.
\end{equation*}

\subsection{Combining identities}
This technique deals with combining the Pythagorean trigonometric identity with the double-angle formula and the identity of $\cos x$ and $\tan (x/2)$.
The double-angle formula used here is $\cos x = \cos^2(x/2) - \sin^2(x/2)$.
The identity of $\cos x$ in terms of $\tan (x/2)$ reads
\begin{equation*}
  \cos x = \frac{1 - \tan^2 (x/2)}{1 + \tan^2 (x/2)}.
\end{equation*}
Employing these identities the integral $J$ now reads
\begin{equation}
  J = \left\{
        \begin{array}{rl}
          2\sqrt{2} \sgn \left[\cos (x/2)\right] \, \sin (x/2) + C, & \hbox{for $+$ sign} \\
         -2\sqrt{2} \sgn \left[\sin (x/2)\right] \, \cos (x/2) + C, & \hbox{for $-$ sign}.
        \end{array}
      \right. \label{combi}
\end{equation}

\subsection{Twice trigonometric substitutions}
Employing the substitution $u = \tan (x/2)$, we have
\begin{equation*}
J = \left\{
    \begin{array}{rl}
    {\displaystyle \int \frac{2\sqrt{2}     \, du}{(1 + u^2)^{3/2}} = \sgn\left(\frac{1}{\sqrt{1 + u^2}}\right) \frac{ 2 \sqrt{2} u}{\sqrt{1 + u^2}} + C}, & \hbox{for $+$ sign} \\
    {\displaystyle \int \frac{2\sqrt{2} |u| \, du}{(1 + u^2)^{3/2}} = \sgn\left(\frac{u}{\sqrt{1 + u^2}}\right) \frac{-2 \sqrt{2}  }{\sqrt{1 + u^2}} + C}, & \hbox{for $-$ sign}.
    \end{array}
    \right.
\end{equation*}
After returning to the initial variable $x$, identical expressions with the ones in~\eqref{combi} will be obtained.

\section{Application: Cardioid} \label{cardioid}
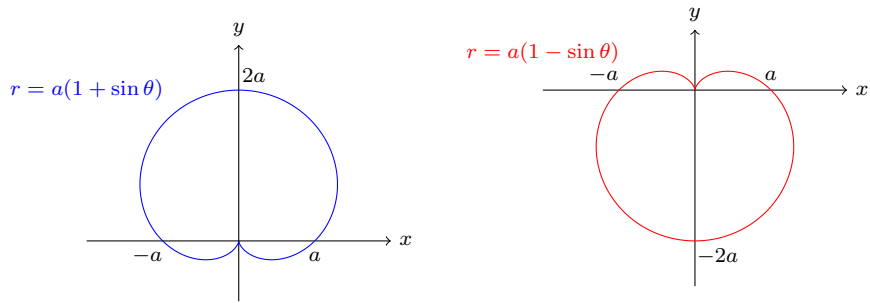
\begin{figure}[h]
\begin{center}
\begin{tikzpicture}[scale=2]
\draw[->] (-1,0) -- (1,0);
\draw[->] (0,-0.4) -- (0,1.3);
\draw node [blue] at (-1,1) {\scriptsize{$r = a(1 + \sin \theta)$}};
\draw node [black] at (0.5,-0.1) {\scriptsize{$a$}};
\draw node [black] at (-0.6,-0.1) {\scriptsize{$-a$}};
\draw node [black] at (0.1,1.1) {\scriptsize{$2a$}};
\draw node [black] at (1.1,0) {\scriptsize{$x$}};
\draw node [black] at (0,1.4) {\scriptsize{$y$}};
\draw[color=blue,domain=0:6.28,samples=200,smooth] plot (xy polar cs:angle=\x r,radius = {.5+.5*sin(\x r)});    

\begin{scope}[shift={(3,1)}]
\draw[->] (-1,0) -- (1,0);
\draw[->] (0,-1.3) -- (0,0.4);
\draw node [red] at (-1,.25) {\scriptsize{$r = a(1 - \sin \theta)$}};
\draw node [black] at (0.5,0.1) {\scriptsize{$a$}};
\draw node [black] at (-0.6,0.1) {\scriptsize{$-a$}};
\draw node [black] at (0.15,-1.1) {\scriptsize{$-2a$}};
\draw node [black] at (1.1,0) {\scriptsize{$x$}};
\draw node [black] at (0,0.5) {\scriptsize{$y$}};
\draw[color=red,domain=0:6.28,samples=200,smooth] plot (xy polar cs:angle=\x r,radius = {.5-.5*sin(\x r)});    
\end{scope}
\end{tikzpicture}
\caption{Sketches of cardioids $r = a(1 + \sin \theta)$ (left) and $r = a(1 - \sin \theta)$ (right), $a > 0$.} \label{plot}
\end{center}
\end{figure}

The integral discussed above appears as calculation of the arc length of a cardioid.
The length of cardioids $r = a(1 \pm \sin \theta)$, $a > 0$ is given by
\begin{equation*}
  L = a \sqrt{2} \int_0^{2\pi} \sqrt{1 \pm \sin \theta}\, d\theta.
\end{equation*}
The sketches of the cardioids are presented in Figure~\ref{plot}.
The properties of the curve have been investigated in a classical paper by Yates more than half a century ago~\cite{yates59}.
The author also compiled a handbook on many kinds of curves, including cardioid, and discussed their properties~\cite{yates47}.
Another approach of calculating an area of cardioid and other shapes of closed curves is presented using the surveyor's method~\cite{braden}.
A road-wheel relationship by rolling a cardioid wheel on an inverted cycloid is discussed in~\cite{hall}.

Cardioid finds various applications in fractals, complex analysis, plant physiology and engineering. 
In fractals, it appears in Douady cauliflower, which is a decoration formed via numerous small cardioids of the Mandelbrot set~\cite{pastor,romera}.
In plant physiology, the seed shape of {\slshape Arabidopsis} (rock cress) can be modelled using cardioid~\cite{cervantes}.
The model based on the comparison of the outline of the seed's longitudinal section with a transformed cardioid,
where the horizontal axis is scaled by a factor equal to the Golden Ratio.
An envelope of rays either reflected or refracted from the surface, known as caustic, from a cup of coffee or milk exhibits the shape of a cardioid~\cite{caustic}. 
In the field of electronics and electrical engineering, a cardioid directional pattern in a microphone provides a relatively wide pick-up zone~\cite{marshall}. 

It is stated but not shown in a Calculus textbook authored by Stewart~\cite{stewart10, stewart12} that one can calculate this integral using the techniques described in this article or by technology, amongst others, are {\sl Integral Calculator}~\cite{calculator}, {\sl Sage}~\cite{sage}, {\sl Symbolab}~\cite{symbolab} and {\sl Wolfram Alpha}~\cite{wolfram}. The author uses the cardioid $r = 1 + \sin \theta$ as an example, as shown in Figure~\ref{example} mentioned earlier in the introduction of this article.
In general, evaluating a definite integral involving an absolute value, one must find the zeros of the function in the absolute value and divide the range of integration into pieces by toggling the sign within each of the intervals.

\subsection{Rationalizing numerator}
\begin{figure}[h]
\begin{center}
\begin{tikzpicture}
\draw[->] (-1,0) -- (7,0);
\draw[->] (0,-0.4) -- (0,1.3);
\draw [color=red] plot[domain=0:2*pi] (\x,{abs(cos(\x r))});
\draw node [black] at (0,1.5) {\scriptsize{$|\cos \theta|$}};
\draw node [black] at (7.1,0) {\scriptsize{$\theta$}};
\draw node [black] at (1.5,-0.2) {\scriptsize{$\pi/2$}};
\draw node [black] at (4.7,-0.2) {\scriptsize{$3\pi/2$}};
\draw node [black] at (-0.2,1) {\scriptsize{$1$}};

\begin{scope}[shift={(0,-3)}]
\draw[->] (-1,0) -- (7,0);
\draw[->] (0,-0.4) -- (0,1.3);
\draw [color=blue] plot[domain=0:2*pi] (\x,{abs(cos(\x/2 r) + sin(\x/2 r))});
\draw node [black] at (0.5,1.7) {\scriptsize{$|\cos \theta/2 + \sin \theta/2|$}};
\draw node [black] at (7.1,0) {\scriptsize{$\theta$}};
\draw node [black] at (4.7,-0.2) {\scriptsize{$3\pi/2$}};
\draw node [black] at (-0.2,1) {\scriptsize{$1$}};
\end{scope}

\begin{scope}[shift={(0,-6)}]
\draw[->] (-1,0) -- (7,0);
\draw[->] (0,-0.4) -- (0,1.3);
\draw [color=blue] plot[domain=0:2*pi] (\x,{abs(cos(\x/2 r) - sin(\x/2 r))});
\draw node [black] at (0.5,1.7) {\scriptsize{$|\cos \theta/2 - \sin \theta/2|$}};
\draw node [black] at (7.1,0) {\scriptsize{$\theta$}};
\draw node [black] at (1.5,-0.2) {\scriptsize{$\pi/2$}};
\draw node [black] at (-0.2,1) {\scriptsize{$1$}};
\end{scope}

\end{tikzpicture}
\caption{Plots of $|\cos \theta|$ (top panel), $\left|\cos (\theta/2) + \sin (\theta/2) \right|$ (middle panel) 
and $\left|\cos (\theta/2) - \sin (\theta/2) \right|$ (bottom panel) for $0 \leq \theta \leq 2\pi$ with the indicated zeros.} \label{zeros}
\end{center}
\end{figure}
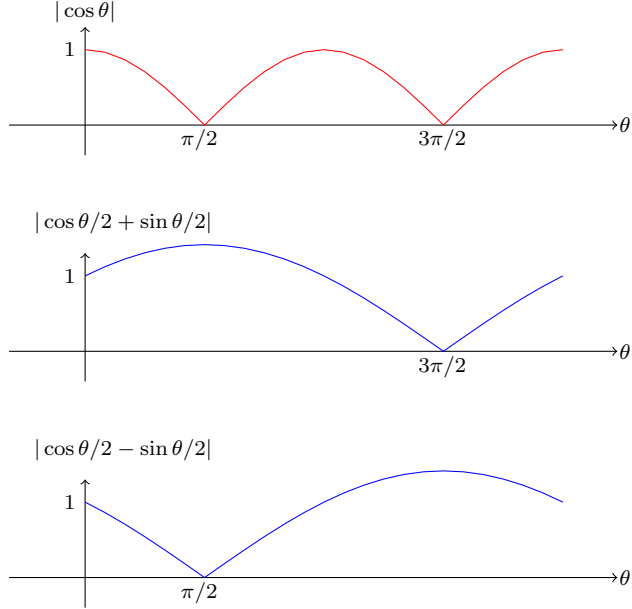

Since $\cos \theta \geq 0$ for $0 \leq \theta \leq \pi/2$ and $3\pi/2 \leq \theta \leq 2\pi$ and $\cos \theta < 0$ for $\pi/2 < \theta < 3\pi/2$,
we need to split the integral into three intervals. 
See the top panel of Figure~\ref{zeros}. Thus, using the result from Subsection~\ref{ratnum}, the length of the cardioids $r = a(1 \pm \sin \theta)$ is given by
\begin{eqnarray*}
L &=& a \sqrt{2} \int_{0}^{2\pi} \frac{|\cos \theta|}{\sqrt{1 \mp \sin \theta}} \, d\theta \\
  &=& a \sqrt{2} \left( \int_0^{\pi/2} \frac{\cos \theta}{\sqrt{1 \mp \sin \theta}} \, d\theta - \int_{\pi/2}^{3\pi/2} \frac{\cos \theta}{\sqrt{1 \mp \sin \theta}} \, d\theta
        + \int_{3\pi/2}^{2\pi} \frac{\cos \theta}{\sqrt{1 \mp \sin \theta}} \, d\theta \right) \\
  &=& 2 a \sqrt{2} \left(\left. \mp \sqrt{1 \mp \sin \theta} \right|_{0}^{\pi/2} \pm \left.\sqrt{1 \mp \sin \theta} \right|_{\pi/2}^{3\pi/2} \mp \left. \sqrt{1 \mp \sin \theta} \right|_{3\pi/2}^{2\pi} \right) 
   =  8a.
\end{eqnarray*}

\subsection{Twice trigonometric substitutions}
We know that (see the middle panel of Figure~\ref{zeros})
\begin{equation*}
  \cos \frac{\theta}{2} + \sin \frac{\theta}{2} = \sqrt{2} \cos \left(\frac{\theta}{2} - \frac{\pi}{4} \right)
  \left\{
    \begin{array}{ll}
      \geq 0, & \hbox{for \; \;} 0 \leq \theta \leq 3\pi/2 \\
      < 0   , & \hbox{for \; \;} 3\pi/2 < \theta \leq 2\pi
    \end{array}
  \right.
\end{equation*}
Thus, implementing this method, the length of the cardioid $r = a(1 + \sin \theta)$ reads
\begin{eqnarray*}
L &=& a \sqrt{2} \int_{0}^{2\pi} \left|\cos \frac{\theta}{2} + \sin \frac{\theta}{2} \right| \, d\theta \\
  &=& a \sqrt{2} \left( \int_0^{3\pi/2} \left( \cos \frac{\theta}{2} + \sin \frac{\theta}{2} \right) \, d\theta -
      \int_{3\pi/2}^{2\pi} \left( \cos \frac{\theta}{2} + \sin \frac{\theta}{2} \right) \, d\theta \right) \\
  &=& 2 a \sqrt{2} \left(\left. \sin \frac{\theta}{2} - \cos \frac{\theta}{2} \right|_{0}^{3\pi/2} -
       \left. \left( \sin \frac{\theta}{2} - \cos \frac{\theta}{2} \right) \right|_{3\pi/2}^{2\pi} \right) 
   =  8a.
\end{eqnarray*}
Similarly, splitting the integral at $\theta = \pi/2$, we also obtain the length $L = 8a$ corresponding to the cardioid $r = a(1 - \sin \theta)$.
See the bottom panel of Figure~\ref{zeros} to observe that the zero of $\cos (\theta/2) - \sin (\theta/2)$ for $0 \leq \theta \leq 2 \pi$ is at $\pi/2$.

\subsection{Variable shift}
These integrals involve the absolute value functions $|\cos (y/2)|$ and $|\sin (y/2)|$, for which in the original variable $\theta$,
both functions are non-negative for $0 \leq \theta \leq 3\pi/2$ and negative for $3\pi/2 < \theta < 2\pi$.
Thus, the length of the cardioid $r = a(1 + \sin \theta)$ reads
\begin{eqnarray*}
  L &=& 2 a \int_0^{2\pi} \left| \sin \left(\frac{\theta}{2} + \frac{\pi}{4} \right) \right| \, d\theta \\
    &=& 2 a \left(\int_0^{3\pi/2}      \sin \left(\frac{\theta}{2} + \frac{\pi}{4} \right) \, d\theta -
                \int_{3\pi/2}^{2\pi} \sin \left(\frac{\theta}{2} + \frac{\pi}{4} \right) \, d\theta\right) \\
    &=& 4 a \left(- \left. \cos \left(\frac{\theta}{2} + \frac{\pi}{4} \right) \right|_{0}^{3\pi/2} + \left. \cos \left(\frac{\theta}{2} + \frac{\pi}{4} \right) \right|_{3\pi/2}^{2\pi}  \right) 
     = 8a
\end{eqnarray*}
or
\begin{eqnarray*}
  L &=& 2 a\int_0^{2\pi} \left| \cos \left(\frac{\pi}{4} - \frac{\theta}{2} \right) \right| \, d\theta \\
    &=& 2 a \left(\int_0^{3\pi/2}      \cos \left(\frac{\pi}{4} - \frac{\theta}{2} \right) \, d\theta -
                \int_{3\pi/2}^{2\pi} \cos \left(\frac{\pi}{4} - \frac{\theta}{2} \right) \, d\theta\right) \\
    &=& 4 a\left(\left. -\sin \left(\frac{\pi}{4} - \frac{\theta}{2} \right) \right|_{0}^{3\pi/2} + \left. \sin \left(\frac{\pi}{4} - \frac{\theta}{2} \right) \right|_{3\pi/2}^{2\pi}  \right) 
    = 8a.
\end{eqnarray*}
Employing a similar technique, identical result of $L = 8a$ is also obtained for the corresponding cardioid $r = a(1 - \sin \theta)$.

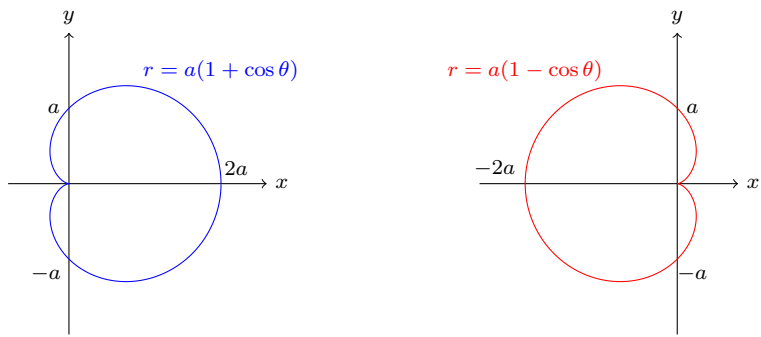
\begin{figure}[h]
\begin{center}
\begin{tikzpicture}[scale=2]
\draw[->] (-0.4,0) -- (1.3,0);
\draw[->] (0,-1) -- (0,1);
\draw node [blue] at (1,0.75) {\scriptsize{$r = a(1 + \cos \theta)$}};
\draw node [black] at (-0.1,0.5) {\scriptsize{$a$}};
\draw node [black] at (-0.15,-0.6) {\scriptsize{$-a$}};
\draw node [black] at (1.1,0.1) {\scriptsize{$2a$}};
\draw node [black] at (0,1.1) {\scriptsize{$y$}};
\draw node [black] at (1.4,0) {\scriptsize{$x$}};
\draw[color=blue,domain=0:6.28,samples=200,smooth] plot (xy polar cs:angle=\x r,radius = {.5+.5*cos(\x r)});    

\begin{scope}[shift={(4,0)}]
\draw[->] (0,-1) -- (0,1);
\draw[->] (-1.3,0) -- (0.4,0);
\draw node [red] at (-1,0.75) {\scriptsize{$r = a(1 - \cos \theta)$}};
\draw node [black] at (0.1,0.5) {\scriptsize{$a$}};
\draw node [black] at (0.1,-0.6) {\scriptsize{$-a$}};
\draw node [black] at (-1.2,0.1) {\scriptsize{$-2a$}};
\draw node [black] at (0.5,0) {\scriptsize{$x$}};
\draw node [black] at (0,1.1) {\scriptsize{$y$}};
\draw[color=red,domain=0:6.28,samples=200,smooth] plot (xy polar cs:angle=\x r,radius = {.5-.5*cos(\x r)});    
\end{scope}
\end{tikzpicture}
\caption{Sketches of cardioids $r = a(1 + \cos \theta)$ (left) and $r = a(1 - \cos \theta)$ (right), $a > 0$.} \label{plot2}
\end{center}
\end{figure}

An expression $r = a(1 \pm \cos \theta)$, $a > 0$ produces cardioids too.
When comparing this expression with the one with sine term, the effect is a 90-degree rotation, either clockwise (for the same sign) or counterclockwise (for the opposite sign), of the corresponding cardioids with the sine term.
The sketch of the corresponding cardioid is presented in Figure~\ref{plot2}.
The length of cardioids $r = a(1 \pm \cos \theta)$, $a > 0$ is given by
\begin{equation*}
  L = a \sqrt{2} \int_0^{2\pi} \sqrt{1 \pm \cos \theta}\, d\theta.
\end{equation*}
Using similar techniques discussed in Section~\ref{cosine}, one can find that the length of these cardioids is also~$8a$.

A number of Calculus textbooks use this type of cardioid as an example for calculating its length. 
For instance, Anton {et al.}~\cite{anton} uses the cardioid $r = 1 + \cos \theta$. After some manipulations, one needs to integrate $|\cos (\theta/2)|$ from $\theta = 0$ to $\theta = 2\pi$. Although general readers will attempt to split the boundary integrations at $\theta = \pi$, the authors explain that since the cardioid is symmetry about the polar axis, the integral from $\theta = \pi$ to $\theta = 2\pi$ is equal to the one from $\theta = 0$ to $\theta = \pi$. Thus, the integral can be calculated by twice integrating from $\theta = 0$ to $\theta = \pi$ of the positive integrand $\cos (\theta/2)$ (without the absolute value). Calculus' Thomas textbook~\cite{thomas} adopts the cardioid $r = 1 - \cos \theta$. The integrand reduces to $|\sin (\theta/2)|$. Fortunately, $\sin (\theta/2) \geq 0$ for $0 \leq \theta \leq 2\pi$ and thus by removing the absolute value and evaluating the integral, one can quickly obtain the length of the cardioid.

\section{Conclusion and Remark}
This article presents the integral with radical sine and cosine functions where its application appears in the length of a cardioid. 
It turns out that several techniques of integration exist to solve the problem, which is interesting from the perspective of teaching and learning mathematics. 
Despite the current trend of using CAS, the collection of integration techniques presented in this article is a valuable pedagogical tool.
To the best of our knowledge, this is the first time such a compilation for this particular type of integrands is presented.
We are convinced that this article contains useful educational contents that will be beneficial for both instructors and students alike.
We also consider our contribution as a complement to existing Calculus textbooks which discuss a topic on calculating the length of a polar curve, particularly cardioid.

\begin{acknowledgment}
{\small \normalfont The authors wish to thank Dr. Ulrich Norbisrath (Faculty of Computer Science, Communication and Media, University of Applied Sciences Upper Austria), 
Dr. Richard J. Mathar
(Max-Planck Institut f\"{u}r Astronomie, Heidelberg, Germany), Professor Victor Hugo Moll (Tulane University, New Orleans, Louisiana, USA),
Professor Chris Sangwin (Mathematics Education Centre, Loughborough University and School of Mathematics, The University of Edinburgh, UK),
the anonymous reviewers whose comments and remarks helped the improvement of this article,
Murat Yessenov (Class 2017 of Physics major, SST, NU) 
and other students from SST and SHSS (School of Humanities and Social Sciences) 
who enrolled in Section~3 of MATH-162 Calculus~2 during Spring~2014 at NU, Astana, Kazakhstan. 
\par}
\end{acknowledgment}

\end{document}